\documentclass[a4paper, 12pt, leqno]{amsart}
\usepackage{mathrsfs}
\usepackage{amsmath}
\usepackage{amscd}
\usepackage{amssymb}
\usepackage{amsthm}
\usepackage{geometry}
\usepackage{tikz}

\title[Based Ring of Two-sided Cells ]
{The Based Rings of Two-sided cells in an Affine Weyl group of type $\tilde B_3$, II}

\author[Y. Qiu and N. Xi]{Yannan Qiu$^{*}$ and Nanhua XI$^{\dagger}$}
\address{$^{*}$
School of Mathematical Sciences\\
Zhejiang University \\
Zhejiang 310058, China \\
China} \email{qiuyannan@zju.edu.cn}
\thanks{Y. Qiu was partially supported by National Natural Science Foundation of China, No. 12171030.}
\address{$^{\dagger}$
Academy of Mathematics and Systems Science\\
Chinese Academy of Sciences\\
Beijing 100190, China\\
and\\
School of Mathematical Sciences\\
University of Chinese Academy of Sciences\\
Chinese Academy of Sciences\\
Beijing 100049, China} \email{nanhua@math.ac.cn}
\thanks{N. Xi was partially supported by National Key R\&D Program of China, No. 2020YFA0712600, and by National Natural Science Foundation of China, No. 11688101.}
\dedicatory{Dedicated to George Lusztig with greatest respect.}

\begin{document}
\baselineskip=18pt
\begin{abstract}
We compute the based rings of two-sided cells corresponding to the unipotent classes
 in $Sp_6(\mathbb C)$ with Jordan blocks (33), (411), (222), respectively. The results for the first two two-sided cells also verify Lusztig's conjecture on the structure of the based rings of two-sided cells of an affine Weyl group. The result for the last two-sided cell partially suggests a modification of Lusztig's conjecture on the structure of the based rings of two-sided cells of an affine Weyl group.
\end{abstract}

\maketitle




We are concerned with the based rings of two-sided cells in an affine Weyl group of type $\tilde B_3$. In a previous paper we discussed the  based ring   of the two-sided cell  corresponding to the nilpotent
element in $Sp_6(\mathbb C)$ with 3 equal Jordan blocks and showed that Lusztig's conjecture on the structure of the based rings of the two-sided cells of an affine Weyl {group}  needs modification (see section 4 in [QX]). In this paper we compute the based rings of two-sided cells corresponding to the unipotent classes
 in $Sp_6(\mathbb C)$ with Jordan blocks (411), (33),  (222),  respectively. The results for the first two two-sided cells  also verify Lusztig's conjecture on the structure of the based rings of two-sided cells of an affine Weyl group. The result for the last two-sided cell partially suggests a modification of Lusztig's conjecture on the structure of the based rings of two-sided cells of an affine Weyl group. For the first two two-sided cells, the validity of Lusztig's conjecture on the based rings is already included in the main theorem in [BO]. Here we construct the bijection in Lusztig's conjecture explicitly so that the results in this paper can be used for computing certain irreducible representations of affine Hecke algebras of type $\tilde B_3$. In section 5 we  give a description for the based ring of the two-sided cell  corresponding to the nilpotent
element in $Sp_6(\mathbb C)$ with  Jordan blocks (222), which  can also be used to compute certain irreducible representations of affine Hecke algebras of type $\tilde B_3$.

The contents of the paper are as follows. Section 1 is devoted to preliminaries, which {include} some basic facts on (extended) affine Weyl groups and their Hecke algebras and {formulation} of  Lusztig's conjecture on the structure of the based ring of a two-sided cell in an affine Weyl group. In section 2 we recall some results on cells of the (extended) affine Weyl group of type $\tilde B_3$, which are due to J. Du.
Sections 3, 4, 5 are devoted to discussing based rings of two-sided cells corresponding to the unipotent classes
 in $Sp_6(\mathbb C)$ with Jordan blocks (411), (33),  (222), respectively.

\def\Cal{\mathcal}
\def\bold{\mathbf}
\def\ca{\mathcal A}
\def\cdz{\mathcal D_0}
\def\cd{\mathcal D}
\def\cdo{\mathcal D_1}
\def\bold{\mathbf}
\def\l{\lambda}
\def\le{\leq}

\def\ll{\underset {L}{\leq}}
\def\rl{\underset {R}{\leq}}
\def\lr{\rl}
\def\lrl{\underset {LR}{\leq}}
\def\llr{\lrl}
\def\el{\underset {L}{\sim}}
\def\er{\underset {R}{\sim}}
\def\elr{\underset {LR}{\sim}}
\def\ds{\displaystyle\sum}

\section{Affine Weyl groups and their Hecke algebras}

In this section we fix some notations and refer to [KL, L1, L2, L3, QX] for more details.

\medskip

 \noindent{\bf 1.1. Extended affine Weyl groups and their Hecke algebras} \
 Let $G$ be a connected reductive algebraic group over
the field {$\mathbb C$} of complex numbers. Let $W_0$ be the Weyl group of $G$ and $W$ the affine Weyl group attached to $G$.   The set of simple reflections of $W$ is denoted by $S$. We shall
denote the length function of $W$ by $l$ and use $\leq$ for the
Bruhat order on $W$. We {also} often write $y<w$ or $w>y$ if $y\le w$ and $y\ne w$.

Let $H$ be the Hecke algebra  of $(W,S)$ over $\Cal A=\mathbb
Z[q^{\frac 12},q^{-\frac 12}]$ $(q$ an indeterminate) with
parameter $q$. Let  $\{T_w\}_{w\in W}$ be its standard basis. Then we have $(T_s-q)(T_s+1)=0$ and $T_wT_u=T_{wu}$ if $l(wu)=l(w)+l(u)$.
Let
$C_w=q^{-\frac {l(w)}2}\sum_{y\le w}P_{y,w}T_y,\ w\in W$ be the
Kazhdan-Lusztig basis of $H$, where $P_{y,w}$ are the
Kazhdan-Lusztig polynomials. The degree of $P_{y,w}$ is less than
or equal to $\frac12(l(w)-l(y)-1)$ if $y< w$ and  $P_{w,w}=1$. Convention: set $P_{y,w}=0$ if $y\not\le w$.

  If $y< w$, we write
$P_{y,w}= \mu(y,w)q^{\frac12(l(w)-l(y)-1)}
+ \text{lower degree terms.}$
We shall write $y\prec w$ if $\mu(y,w)\ne 0$.  We have

\medskip
\def\vp{\varphi}
\def\st{\stackrel}
\def\sc{\scriptstyle}

(a) Let $y\le w$. Assume that $sw\le w$ for some $s\in S$. Then
\begin{alignat*}{2} P_{y,w}&=P_{sy,w},\ \text{if}\ sy>y;\\
P_{y,w}&=q^{1-c}P_{sy,sw}+q^cP_{y,sw}-\sum_{\st {\st {z\in W}{y\le z\prec sw}}{sz<z}}\mu(z,sw)q^{\frac{l(w)-l(z)}2}P_{y,z},\end{alignat*}
where $c=1$ if $sy< y$ and $c=0$ if $sy>y$.

\medskip

(b) Let $y\le w$. Assume that $ws\le w$ for some $s\in S$. Then
\begin{alignat*}{2}P_{y,w}&=P_{ys,w},\ \text{if}\ ys>y;\\\
P_{y,w}&=q^{1-c}P_{ys,ws}+q^cP_{y,ws}-\sum_{\st {\st {z\in W}{y\le z\prec ws}}{zs<z}}\mu(z,ws)q^{\frac{l(w)-l(z)}2}P_{y,z},\end{alignat*}
where $c=1$ if $ys< y$ and $c=0$ if $ys>y$.

\medskip

From the two formulas above one gets (see [KL])

\medskip

(c) Let $y,w\in W$ and $s\in S$ be such that $y<w,\ sw<w,$ and $ sy>y$. Then $y\prec w$ if and only if $w=sy$. Moreover this implies that $\mu(y,w)=1$.

\medskip

(d) Let $y,w\in W$ and $s\in S$ be such that $y<w,\ ws<w$, and $ys>y$. Then $y\prec w$ if and only if $w=ys$. Moreover this implies that $\mu(y,w)=1$.

\medskip

The following  formulas for computing $C_w$ (see [KL]) will be used in sections 3, 4, 5.

\medskip

(e) Let $w\in W$ and $s\in S${. Then}
\begin{align} C_sC_w=\begin{cases}\displaystyle (q^{\frac12}+q^{-\frac12})C_w,\quad &\text{if\ }sw<w,\\
\displaystyle C_{sw}+\sum_{\st  {z\prec w}{sz<z}}\mu(z,w)C_z,\quad&\text{if\ }sw\ge w.\end{cases}\end{align}
\begin{align} C_wC_s=\begin{cases}\displaystyle (q^{\frac12}+q^{-\frac12})C_w,\quad &\text{if\ }ws<w,\\
\displaystyle C_{ws}+\sum_{\st  {z\prec w}{zs<z}}\mu(z,w)C_z,\quad&\text{if\ }ws\ge w.\end{cases}\end{align}

\medskip

\noindent{\bf 1.2. Cells of affine Weyl groups}  We refer to [KL] for definition of left cells, right cells and two-sided cells of $W$.

For $h,\, h'\in H$ and $x\in W$, write
$$ hC_x =\sum_{y\in W}a_yC_y,\quad
 C_xh =\sum_{y\in W}b_yC_y,\quad
  hC_xh'=\sum_{y\in W}c_yC_y,\quad a_y,b_y, c_y\in \mathcal A.$$
  Define $y\ll x$ if $a_y\ne 0$ for some $h\in H$, $y\rl x$ if $b_y\ne 0$ for some $h\in H$, {and} $y\lrl x$ if $c_y\ne 0$ for some $h,h'\in H$.

  We write $x\el y$ if $x\ll y\ll x$,  $x\er y$ if $x\rl y\rl x$, and $x\elr y$ if $x\lrl y\lrl x$. Then $\el,\ \er,\ \elr$ are equivalence relations on $W$. The equivalence classes are called left cells, right cells, and two-sided cells of $W$ respectively. Note that if $\Gamma$ is a left cell of $W$, then $\Gamma^{-1}=\{ w^{-1}\,|\, w\in\Gamma\}$ is a right cell.

\medskip

For $w\in W$, set $R(w)=\{s\in S\,|\, ws\le w\}$ and $L(w)=\{s\in S\,|\, sw\le w\}.$ Then we have (see [KL])

\medskip

(a) $R(w)\subset R(u)$ if $u\ll w$ and $L(w)\subset L(u)$ if $u\rl w.$ In particular, $R(w)= R(u)$ if $u\el w$ and $L(w)= L(u)$ if $u\er w.$

\medskip

\noindent{\bf 1.3.   $*$-operations}\ \ The $*$-operation introduced in [KL] and   generalized in [L1] is a useful tool in the theory of cells of Coxeter groups.

Let $s,t$ be simple reflections in $S$ and assume that $st$ has order  $m\ge 3$. Let $w\in W$ be such that $sw\ge w,\ tw\ge w$. The $m-1$ elements $sw,\ tsw,\ stsw,\ ...,$ is called a left string (with respect to $\{s,t\})$, and the $m-1$ elements $tw,\ stw,\ tstw,\ ..., $ is also called a left string (with respect to $\{s,t\})$. Similarly we define right strings (with respect to $\{s,t\}$). Then (see [L1])

\medskip

(a) A left string in $W$ is contained in a left cell of $W$ and a right string in $W$ is contained in a right cell of $W$.

\medskip
 Assume that $x$ is in a left (resp. right ) string (with respect to $\{s,t\})$ of length $m-1$ and is the $i$th element of the left (resp. right) string,  define ${}^*x$ (resp. $x^*$) to be the $(m-i)$th element of the string, where $*=\{s,t\}$. The following result is proved in [X2].

 \medskip

 (b) Let $x$ be in $W$ such that $x$ is in a left string with respect to $*=\{s,t\}$ and is also in a right string with respect to $\star=\{s',t'\}$. Then ${}^*x$ is in a right string with respect to $\{s',t'\}$ and $x^\star$ is in a left string with respect to $\{s,t\}$. Moreover ${}^*(x^\star)=({}^*x)^\star$. We shall write ${}^*x^\star$
 for ${}^*(x^\star)=({}^*x)^\star$.

\medskip

The following result is due to Lusztig [L1].

\medskip

(c) Let $\Gamma$ be a left cell of $W$ and an element $x\in\Gamma$ is in a right string $\sigma_x$ with respect to $*=\{s,t\}$. Then any element $w\in\Gamma$ is in a right string $\sigma_w$ with respect to $*=\{s,t\}$. Moreover $\Gamma^*=\{w^*\,|\, w\in \Gamma\}$ is a left cell of $W$ and $\Omega=\displaystyle \left(\cup_{w\in\Gamma}\sigma_w\right)-\Gamma$ is a union of at most $m-2$ left cells.

\medskip

Following Lusztig [L1] we set $\tilde\mu(y,w)=\mu(y,w)$ if $y< w$ and $\tilde\mu(y,w)=\mu(w,y)$ if $w<y$. For convenience we also set $\tilde\mu(y,w)=0$ if $y\nless w$ and $w\nless y$. Assume that $x_1,x_2,...,x_{m-1}$  and $y_1,y_2,...,y_{m-1}$ are two left strings with respect to $*=\{s,t\}$. Define
$$a_{ij}=\begin{cases}\tilde\mu(x_i,y_j),\quad &\text{if\ } \{s,t\}\cap L(x_i)=\{s,t\}\cap L(y_j),\\
0,\quad &\text{otherwise}.\end{cases}$$
Lusztig proved the following identities (see Subsection 10.4 in [L1]).

\medskip

(d) If $m=3$, then $a_{11}=a_{22}$ and $a_{12}=a_{21}$.

\medskip

(e) If $m=4$, then
\begin{align} a_{11}=a_{33},\ a_{13}=a_{31},\ a_{22}=a_{11}+a_{13},\ a_{12}=a_{21}=a_{23}=a_{32}.\end{align}

\medskip

\noindent{\bf 1.4. Lusztig's $a$-function}\quad
For $x,y\in W$, write $$C_xC_y=\sum_{z\in
W}h_{x,y,z}C_z,\qquad h_{x,y,z}\in \mathcal A
 =\mathbb Z[q^{\frac 12},q^{-\frac 12}].$$
 Following Lusztig ([L1]), we define
 $$a(z)={\rm min}\{i\in\bold N\ |\ q^{-\frac i2}h_{x,y,z}\in\mathbb Z[q^{-\frac
 12}]{\rm\ for \ all\ }x,y\in W\}.$$ If for any $i$,
 $q^{-\frac i2}h_{x,y,z}\not\in\mathbb Z[q^{-\frac
 12}]{\rm\ for \ some\ }x,y\in W$, we set $a(z)=\infty.$
 The following properties are proved in [L1].

 \medskip

 (a)  We have $a(w)\le
l(w_0)$ for any $w\in W$, where $w_0$ is the longest element in the Weyl group $W_0$.

\medskip

(b) $a(x)\ge a(y)$ if $x\lrl y$. In particular, $a(x)=a(y)$ if $x\elr y$.

\medskip

(c) $x\el y$ (resp. $x\er y,\ x\elr y$) if $a(x)=a(y)$ and $x\ll y$ (resp. $x\lr y,\ x\llr y$).

(d) If $h_{x,y,z}\ne 0$, then $z\lr x$ and $z\ll y$. In particular, $a(z)>a(x)$ if $z\not \er x$, and $a(z)>a(y)$ if $z\not \el y$.

Following Lusztig,
 we define  $\gamma_{x,y,z}$ by the following formula,
 $$h_{x,y,z}=\gamma_{x,y,z}q^{\frac {a(z)}2}+
 {\rm\ lower\ degree\ terms}.$$
Springer showed that $l(z)\ge a(z)$ (see [L2]). Let $\delta(z)$ be
the
 degree of $P_{e,z}$, where $e$ is the neutral element of $W$.
 Then actually one has $l(z)-a(z)-2\delta(z)\ge 0$ (see [L2]). Set
 $$\cd =\{z\in W\ |\ l(z)-a(z)-2\delta(z)=0\}.$$

 The elements of $\cd$ are involutions, called distinguished involutions of
 $(W,S)$ (see [L2]). The following properties are due to Lusztig [L2] except the  (j) (which is trivial) and (k) (proved in [X2]).

 (e) $\gamma_{x,y,z}\ne 0\Longrightarrow x\el y^{-1},\ y\el z,\ x\er z.$

 (f) $x\el y^{-1}$ if and only if $\gamma_{x,y,z}\ne 0$ for some $z\in W$.

 (g) $\gamma_{x,y,z}=\gamma_{y,z^{-1},x^{-1}}=\gamma_{z^{-1}, x,y^{-1}}$.

 (h) $\gamma_{x,d,x}=\gamma_{d,x^{-1},x^{-1}}=\gamma_{x^{-1},x,d}=1$ if $x\el d$ and $d$ is a distinguished involution.

 (i) $\gamma_{x,y,z}=\gamma_{y^{-1},x^{-1},z^{-1}}.$

 (j) If $\omega,\tau\in W$ has length 0, then
 $$\gamma_{\omega x,y\tau,\omega z\tau}=\gamma_{x,y,z},\ \
 \gamma_{ x\omega,\tau y, z}=\gamma_{x,\omega\tau y,z}.$$

 (k) Let $x,y,z\in W$ be such that (1) $x$ is in a left string with respect to $*=\{s,t\}$ and also in a right string with respect to $\#=\{s',t'\}$, (2) $y$ is in a left string with respect to $\#=\{s',t'\}$ and also in a right string with respect to $\star=\{s'',t''\}$, (3) $z$ is in a left string with respect to $*=\{s,t\}$ and also  in a right string with respect to $\star=\{s'',t''\}$. Then
 $$\gamma_{x,y,z}=\gamma_{{}^*x^\#,{}^\#y^\star,{}^*z^\star}.$$

\medskip
\def\tt{\tilde T}

For $w\in W$, set $\tilde T_w=q^{-l(w)/2}T_w$. For $x,y\in W$, write
$$\tt_x\tt_y=\sum_{z\in W}f_{x,y,z}\tt_z,\qquad f_{x,y,z}\in\mathcal A=\mathbb Z[q^{\frac12},q^{-\frac12}].$$

(l) If $x,y,w$ are in a two-sided cell of $W$, $f_{x,y,w}=\lambda q^{\frac{a(w)}2}+$ lower degree terms and as Laurent polynomials in $q^{\frac12}$, deg$f_{x,y,z}\le a(w)$ for all $z\in W$, then
$$\gamma_{x,y,w}=\lambda.$$

 (m) Each left cell (resp. each right cell) of $W$ contains a unique distinguished involution.

 (n) Each two-sided cell of $W$ contains only finitely many left cells.

(o) Let $I$ be a subset of $S$ such that the subgroup $W_I$ of $W$ generated by $I$ is finite. Then the longest element $w_I$ is a distinguished involution.

\medskip

Let $d$ be a distinguished involution in $W$.

(p) For any $\omega\in\Omega$, the element $\omega d\omega^{-1}$ is a distinguished involution.

(q) Suppose $s,t\in S$ and $st$ has order 3. Then $d\in D_{L}(s,t)$ if and only if $d\in D_{R}(s,t)$. If $d\in D_{L}(s,t)$, then ${}^*d^*$ is a distinguished involution.

\medskip

{\bf 1.5.} Assume $s, t\in S$ and $st$ has order 4. Let $w, u, v$ be in $W$ such that $l(ststw)=l(w)+4$ and $l(ststv)=l(v)+4$. We have (see [X2, 1.6.3])
\begin{itemize}
\item[(a)] $\gamma_{tsw, u, tv}= \gamma_{sw, u, stv},$
\item[(b)] $\gamma_{tsw, u, tsv}= \gamma_{sw, u, sv}+ \gamma_{sw, u, stsv},$
\item[(c)] $\gamma_{tsw, u, tstv}= \gamma_{sw, u, stv},$
\item[(d)]$\gamma_{tstw, u, tv}+\gamma_{tw, u, tv}=\gamma_{stw, u, stv},$
\item[(e)]$\gamma_{tstw, u, tsv}=\gamma_{stw, u, stsv},$
\item[(f)]$\gamma_{tstw, u, tstv}+\gamma_{tw, u, tstv}=\gamma_{stw, u, stv}.$
\end{itemize}

\medskip

Assume $s, t\in S$ and $st$ has order 4. Let $w, u, v$ be in $W$ such that $l(ustst)=l(u)+4$ and $l(vstst)=l(v)+4$. We have (loc.cit)
\begin{itemize}
\item[(a')] $\gamma_{w, ut, vst}= \gamma_{w, uts, vs},$
\item[(b')] $\gamma_{w, ust, vst}= \gamma_{w, us, vs}+ \gamma_{w, usts, vs},$
\item[(c')] $\gamma_{w, utst, vst}= \gamma_{w, uts, vs},$
\item[(d')]$\gamma_{w, ut, vtst}+\gamma_{w, ut, vt}=\gamma_{w, uts, vts},$
\item[(e')]$\gamma_{w, ust, vtst}=\gamma_{w, usts, vts},$
\item[(f')]$\gamma_{w, utst, vtst}+\gamma_{w, utst, vt}=\gamma_{w, uts, vts}.$
\end{itemize}

\medskip

\def\ll{\underset {L}{\leq}}
\def\rl{\underset {R}{\leq}}

\def\lrl{\underset {LR}{\leq}}
\def\llr{\lrl}
\def\el{\underset {L}{\sim}}
\def\er{\underset {R}{\sim}}
\def\elr{\underset {LR}{\sim}}
\def\ds{\displaystyle\sum}

\def\vp{\varphi}
\def\st{\stackrel}
\def\sc{\scriptstyle}

 \noindent{\bf 1.6. The based ring of a two-sided cell}\quad For each two-sided cell $c$ of $W$, let $J_c$ be the free $\mathbb Z$-module with a basis $t_w,\ w\in c$. Define
  $$t_xt_y=\sum_{z\in c}\gamma_{x,y,z}t_z.$$
  Then $J_c$ is an associative ring with unit $\sum_{d\in\cd\cap c}t_d.$

  The ring $J=\bigoplus_{c}J_c$ is a ring with unit $\sum_{d\in\cd}t_d$. Sometimes $J$ is called asymptotic Hecke algebra since Lusztig established an injective $\ca$-algebra homomorphism
  $$  \phi: H \to J\otimes\ca,\quad
  C_x \mapsto\sum_{\st {\st {d\in\cd}{w\in W}}{w\el d}}h_{x,d,w}t_w.$$

  \medskip

  \noindent{\bf 1.7. Lusztig's conjecture on the structure of $J_c$}\quad In [L3] Lusztig states a conjecture on $J_c$ using equivariant $K$-groups on finite sets.

  Let $G$ be a  connected reductive  group over $\mathbb C$. Lusztig establishes a bijection between the two-sided cells of the extended affine Weyl group $W$ and the unipotent classes of $G$.

  For each two-sided cell $c$ of $W$, let $u$ be a unipotent element in the unipotent class corresponding to $c$ and let $F_c$ be a maximal reductive subgroup of the centralizer of $u$ in $G$.

  {\bf Conjecture} (Lusztig [L3]): Assume that $G$ is {a} simply connected simple algebraic group over $\mathbb C$. Then there exists a finite set $Y$ with an algebraic action of $F_c$ and a  bijection
  $$\pi: c\to \text{{the set of} isomorphism classes of irreducible $F_c$-vector bundles  on}\ Y\times Y.$$
  such that

  (i)  The bijection $\pi$ induces a ring homomorphism
  $$\pi: J_c\to K_{F_c}(Y\times Y),\ \ t_x\mapsto \pi(x).$$

  (ii) $\pi(x^{-1})_{(a,b)}=\pi(x)_{(b,a)}^*$ is the dual representation of $\pi(x)_{(b,a)}.$

  \section{Cells in an extended affine Weyl group of type $\tilde B_3$}

  In this section $G=Sp_6(\mathbb C)$, so that the extended affine Weyl group $W$ attached to $G$ is of type $\tilde B_3$. The left cells and two-sided cells are described by J. Du (see [D]). We recall his results.

  {\bf 2.1. The Coxeter graph of $W$}. As usual, we number the 4 simple reflections $s_0,\ s_1,\ s_2,\ s_3$ in $W$ so that
  \begin{alignat*}{2} &s_0s_1=s_1s_0,\quad s_0s_3=s_3s_0,\quad s_1s_3=s_3s_1,\\
  &(s_0s_2)^3=(s_1s_2)^3=e,\quad (s_2s_3)^4=e,\end{alignat*}
  where $e$ is the neutral element in $W$. The relations among the simple reflections can be read through the following Coxeter graph:

 
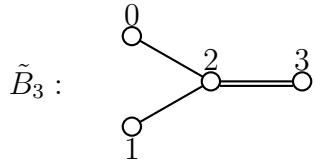
\begin{center}
  \begin{tikzpicture}[scale=.6]
    \draw (-1,0) node[anchor=east]  {$\tilde B_3:$};
    \draw[thick] (2 cm,0) circle (.2 cm) node [above] {$2$};
    \draw[xshift=2 cm,thick] (150:2) circle (.2 cm) node [above] {$0$};
    \draw[xshift=2 cm,thick] (210:2) circle (.2 cm) node [below] {$1$};
    \draw[thick] (4 cm,0) circle (.2 cm) node [above] {$3$};
    \draw[xshift=2 cm,thick] (150:0.2) -- (150:1.8);
    \draw[xshift=2 cm,thick] (210:0.2) -- (210:1.8);
    \draw[thick] (2.2,0) --+ (1.6,0);
    \draw[thick] (2.2,-0.1) --+ (1.6,0);
 \end{tikzpicture}
\end{center}

  There is a unique nontrivial element $\tau $ in $W$ with length 0. We have $\tau^2=e,\ \tau s_0\tau=s_1,\ \tau s_i\tau=s_i$ for $i=2,3.$ Note that $s_1,s_2,s_3$ generate the Weyl group $W_0$ of type $B_3$ and $s_0,s_1,s_2,s_3$ generate an affine Weyl group $W'$ of type $\tilde B_3$. And $W$ is generated by $\tau,\ s_0,s_1,s_2,s_3$.

  \medskip

  {\bf 2.2. Cells in $W$}\quad According to [D], the extended affine Weyl group $W$ attached to  $Sp_6(\mathbb C)$ has 8 two-sided cells:
    $$A,\quad B,\quad C,\quad D,\quad E, \quad F, \quad G,\quad H.$$
   The following table displays some useful information on these two-sided cells.
  \begin{center}
  \doublerulesep 0.4pt \tabcolsep8pt
\begin{tabular}{ccccc}
\hline
& &Number & Size of Jordan blocks& Maximal reductive subgroup\\
& & of left & of the  corresponding & of the centralizer of a unipotent \\
$X$ & {${a(X)}$} & cells in $X$ &unipotent class in $Sp_6(\mathbb C)$&element in the corresp. unipotent class\\
\hline$A$&9&48&\hfill (111111)\ \ \ \ \ \ \ \ \ \ \ \   &$Sp_6(\mathbb C)$\\
$B$&6&24&\hfill (21111)\ \ \ \ \ \ \ \ \ \ \ \  &$Sp_4(\mathbb C)\times\mathbb Z/2\mathbb Z$\\
$C$&4&18&\hfill (2211)\ \ \ \ \ \ \ \ \ \ \ \  &$SL_2(\mathbb C)\times O_2(\mathbb C)$\\
$D$&3&12&\hfill (222)\ \ \ \ \ \ \ \ \ \ \ \  &$O_3(\mathbb C)$\\
$E$&2&8&\hfill (411)\ \ \ \ \ \ \ \ \ \ \ \  &$SL_2(\mathbb C)\times\mathbb Z/2\mathbb Z$\\
$F$&2&6&\hfill (33)\ \ \ \ \ \ \ \ \ \ \ \  &$SL_2(\mathbb C)$\\
$G$&1&4&\hfill (42)\ \ \ \ \ \ \ \ \ \ \ \  &$\mathbb Z/2\mathbb Z\times\mathbb Z/2\mathbb Z$\\
$H$&0&1&\hfill (6)\ \ \ \ \ \ \ \ \ \ \ \  &$\mathbb Z/2\mathbb Z$\\
\hline
\end{tabular}
\end{center}\vspace{2mm}
The notations for two-sided cells in the table are the same as those in [D], which will be replaced by other notations in subsequent sections, otherwise confusion would happen since notations $C, F,G$ are already used for other objects.

In subsequent sections, for a reduced expression  $s_{i_1}s_{i_2}\cdots s_{ik}$ of an element in $W$, we often write $i_1i_2\cdots i_k$ instead of the reduced expression.

{{In}} the rest of the paper, $W$ always stands for the affine Weyl group attached to $Sp_6(\mathbb C)$, $\tau,\ s_i$ are as in Subsection 2.1, and all representations in this paper are rational representations of algebraic groups.

 \section{The based ring of the two-sided cell containing $s_0s_1$}

 {{{\bf 3.1.}} In this section $c$ stands for the two-sided cell of $W$ containing $s_0s_1$. According to [D, Figure I, Theorem 6.4], $c$ has six left cells. We list the six left cells and representative elements in the left cells given in [D, figure I]:
 $$\Gamma_{1},\ 01;\quad \Gamma_2,\ 012;\quad \Gamma_3,\ 0123;\quad \Gamma_4,\ 01232;\quad \Gamma_5,\ 012321;\quad\Gamma_6,\ 012320.$$
  The values of $a$-function on $c$ is 2.

 The corresponding unipotent class in $Sp_6(\mathbb C)$ has Jordan block sizes (411). {Maximal reductive subgroup} of the centralizer of an element in the unipotent class is $F_c=\mathbb Z/2\mathbb Z\times SL_2(\mathbb C)$. Let $\epsilon$ be the nontrivial one dimensional representation of $ \mathbb Z/2\mathbb Z$ and $V(k)$ be an irreducible representation of $SL_2(\mathbb C)$ with highest weight $k$. They can be regarded as irreducible representations of $F_c$ naturally. Up to isomorphism, the irreducible representations of $F_c$ are $V(k),\ \epsilon\otimes V(k),\ k=0,1,2,3,...$. We will denote $ \epsilon\otimes V(k)$ by  $\epsilon  V(k)$.

\def\Rep{\text{Rep\,}}

 Let $x_k=(s_0s_1s_2s_3s_2)^ks_1s_0, \ u_1=e,\ u_2=s_2,\ u_3=s_3s_2,\ u_4=s_2s_3s_2,\ u_5=s_1s_2s_3s_2,\ u_6=s_0s_2s_3s_3.$

 According to [D,  Theorem 6.4], we have

(a) $ c=\{ u_ix_ku_j^{-1},\   u_i \tau x_ku_j^{-1}\,|\, 1\le i,j\le 6, \ k=0,1,2,3,...\}.$

(b) $\Gamma_j=\{ u_ix_ku_j^{-1},\  u_i\tau x_ku_j^{-1}\,|\, 1\le i \le 6, \ k=0,1,2,3,...\},\ j=1,2,3,4,5,6.$

Let $Y=\{1,\ 2, \ ...,\ 6\}$ and let $F_c$ act on $Y$ trivially. Then $K_{F_c}(Y\times Y)$ is isomorphic to the $6\times 6$ matix ring $M_6(\Rep F_c)$, where $\Rep F_c$ is the representation ring of $F_c$. Recall that $F_c=\mathbb Z/2\mathbb Z\times SL_2(\mathbb C)$ in this section.

The main result in this section is the following theorem.

\medskip

{\bf Theorem 3.2.} Let $c$ be the two-sided cell of $W$ (the affine Weyl group $W$ attached to $Sp_6(\mathbb C)$) containing $s_0s_1$. Then the map
$$\pi: c\to M_6(\Rep F_c), \quad u_ix_ku_j^{-1}\mapsto V(k)_{ij},\  u_i\tau x_ku_j^{-1}\mapsto \epsilon V(k)_{ij}\ $$
induces a ring isomorphism
$$\pi: J_c\to M_6(\Rep F_c),\quad t_{u_ix_ku_j^{-1}}\mapsto V(k)_{ij},\quad t_{u_i\tau x_ku_j^{-1}}\mapsto \epsilon V(k)_{ij},$$
where $V(k)_{ij}$ (resp. $\epsilon V(k)_{ij})$ is the matrix in $M_6(\Rep F_c)$ whose entry at $(p,q)$ is $V(k)$ (resp. $\epsilon V(k)$) if $(p,q)=(i,j)$ and is 0 otherwise.

\medskip

{\sl Remark:} The Theorem 4 in [BO] implies that Lusztig's conjecture on the structure of $J_c$ is true. Since under the isomorphism $K_{F_c}(Y\times Y)\simeq M_6(\Rep F)$, irreducible ${F_c}$-vector bundles on $Y\times Y$ correspond to those $V(k)_{ij},\ \epsilon V(k)_{ij}$, hence Theorem 3.2 provides a computable verification for Lusztig conjecture on the structure  of $J_c$.

\medskip

We prove Theorem 3.2 by  establishing three lemmas.

\medskip

{\bf Lemma 3.3.} Let $1\le i,j,m,n\le 6$ and $k,l$ be nonnegative integers{. For} $z_k=x_k$ or $\tau x_k$, {$z_l=x_l$ or $\tau x_l$, and} $z_p=x_p$ or $\tau x_p$, we have

(a) $\gamma_{u_iz_ku_j^{-1},u_mz_lu_n^{-1},z}=0\quad \text{if}\ j\ne m\ \text{or}\ z\ne u_i\tau^a z_pu^{-1}_n,\ a=0,1,\ \text{for some }p;$

(b) $\gamma_{u_iz_ku_j^{-1},u_jz_lu_n^{-1},u_iz_pu_n^{-1}}=\gamma_{z_k,z_l,z_p},\quad\text{for any nonnegative integer } p.$

\medskip

{\bf Proof.} Note that $z_l^{-1}=z_l$. If $ \gamma_{u_iz_ku_j^{-1},u_mz_lu_n^{-1},z}\ne0$, then by 1.4(e) we get $u_iz_ku_j^{-1}\el(u_mz_lu_n^{-1})^{-1}=u_nz_lu_m^{-1}$, $u_iz_ku_j^{-1}\er z,\ u_mz_lu_n^{-1}\el z$. By (b) in Subsection 3.1 we see that the first assertion is true.

Now we prove the second assertion. Let $ *=\{s_1,s_2\},\ \#=\{s_2,s_3\}$ and $\star=\{s_0,s_2\}$. Then

\medskip

(c) $\Gamma_2=\Gamma_1^*,\quad \Gamma_4=\Gamma_2^{\#},\quad \Gamma_5=\Gamma_4^*,\quad \Gamma_6=\Gamma_4^\star.$

\medskip

Applying 1.4 (k) we see that (b) is true if none of $i,j,n$ is 3.

Now assume that $i=3$. By 1.5 (b)  we get
$$\gamma_{u_3z_ku_j^{-1},u_jz_lu_n^{-1},u_3z_pu_n^{-1}}=\gamma_{u_2z_ku_j^{-1},u_jz_lu_n^{-1},u_2z_pu_n^{-1}}
+\gamma_{u_2z_ku_j^{-1},u_jz_lu_n^{-1},u_4z_pu_n^{-1}}.$$
By  Part (a) of the lemma, we have $\gamma_{u_2z_ku_j^{-1},u_jz_lu_n^{-1},u_4z_pu_n^{-1}}=0.$ Then using (c) above and 1.4 (k) we get
$$\gamma_{u_3z_ku_j^{-1},u_jz_lu_n^{-1},u_3z_pu_n^{-1}}=\gamma_{u_2z_ku_j^{-1},u_jz_lu_n^{-1},u_2z_pu_n^{-1}}
=\gamma_{ z_ku_j^{-1},u_jz_lu_n^{-1}, z_pu_n^{-1}}.$$

Similarly, if $n=3$, we have
$$\gamma_{u_iz_ku_j^{-1},u_jz_lu_3^{-1},u_iz_pu_3^{-1}}=
\gamma_{u_iz_ku_j^{-1},u_jz_lu_2^{-1},u_iz_pu_2^{-1}}
=\gamma_{ u_iz_ku_j^{-1},u_jz_l, u_iz_p}.$$

We have showed for any $1\le i,n\le 6$ the following identity holds:
$$\gamma_{u_iz_ku_j^{-1},u_jz_lu_n^{-1},u_iz_pu_n^{-1}}=\gamma_{z_ku_j^{-1},u_j z_l,z_p}.$$

Note that $z_k^{-1}=z_k$. By above identity and  1.4 (g), we get
$$\gamma_{z_ku_j^{-1},u_jz_l,z_p}=\gamma_{u_jz_l,z_p^{-1},u_jz_k^{-1}}=\gamma_{z_l,z_p,z_k}=\gamma_{z_k,z_l,z_p}.$$
Assertion (b) is proved and the lemma is proved.\qed

\medskip

{\bf Lemma 3.4.} For nonnegative integers, and $a,b =0,1$, we have
$$\gamma_{\tau^a x_k,\tau ^b x_l,\tau^{a+b} x_p}=\gamma_{x_k,x_l,x_p}, \quad \gamma_{\tau^a x_k,\tau ^b x_l,\tau^{c} x_p}=0\ \text{if }\tau^c\ne \tau^{a+b}.$$

\medskip

{\bf Proof.} The assertion follows from  1.4 (j).\qed

\medskip

{\bf Lemma 3.5.} For nonnegative integers $k,l$ we have
$$t_{x_k}t_{x_l}=\sum_{0\le p\le\min\{ k,l\}}t_{{x_{k+l-2i}}}.$$

\medskip

{\bf Proof.} If $k=0$ or $l=0$, the identity above is trivial since $x_0$ is {a distinguished involution}.

Now assume that $k=1$ and $l\ge 1$.  Let $\zeta=q^{\frac12}-q^{-\frac12}$. By a simple computation we see
      $$\tilde T_{x_1}\tilde T_{x_l}=\zeta^2(\tilde T_{x_{l+1}}+\tilde T_{x_{l-1}}+\tilde T_{s_0s_1s_3(s_2s_3s_2s_0s_1)^l}+\tilde T_{s_0s_2s_3s_2s_1(s_2s_3s_2s_0s_1)^l})+\text{lower degree terms,}$$
 Since  $a(s_0s_1s_3(s_2s_3s_2s_0s_1)^l)\geq a(s_0s_1s_3)=3,\  a(s_0s_2s_3s_2s_1(s_2s_3s_2s_0s_1)^l)\geq a(s_2s_1s_2)=3$, {we see that} $s_0s_1s_3(s_2s_3s_2s_0s_1)^l$ and $s_0s_2s_3s_2s_1(s_2s_3s_2s_0s_1)^l$ are not in the two-sided cell $c$. By  1.4 (l), we have
 $$t_{x_1}t_{x_l}=t_{x_{l+1}}+t_{x_{l-1}}.$$

 For $k\ge 2$, since $t_{x_k}=t_{x_1}t_{x_{k-1}}-t_{x_{k-2}}$,  we can use induction on $k$ to prove  the lemma. The argument  is completed.\qed

 \medskip

 {\bf Proof of Theorem 3.2.} Combining Lemmas 3.3, 3.4 and 3.5 we see that Theorem 3.2 is true.

 \section{The based ring of the two-sided cell containing $s_1s_3$}

 {{\bf 4.1.}} In this section $c$ stands for the two-sided cell of $W$ containing $s_1s_3$. According to  [D, Figure I, Theorem 6.4], $c$ has eight left cells. We list the eight left cells and representative elements in the left cells given in [D, Figure I]:
 $$\begin{array}{lllllllll}&\Gamma_{1},\ &13;\quad &\Gamma_2,\ &132;\quad &\Gamma_3,\ &1323;\quad &\Gamma_4,\ &1320;\\
 &\Gamma_{5},\ &03;\quad &\Gamma_6,\ &032;\quad &\Gamma_7,\ &0323;\quad &\Gamma_8,\ &0321.\end{array}$$
 The values of $a$-function on $c$ is 2.

 The corresponding unipotent class in $Sp_6(\mathbb C)$ has Jordan block sizes (33). Maximal reductive subgroup of the centralizer of an element in the unipotent class is ${F_c}= SL_2(\mathbb C)$. Let  $V(k)$ be an irreducible representation of ${F_c}=SL_2(\mathbb C)$ with highest weight $k$.  Up to isomorphism, the irreducible representations of ${F_c}$ are $V(k),\ k=0,1,2,3,...$.

 Let $x_k=(\tau s_0s_3s_2)^ks_1s_3, \ u_1=e,\ u_2=s_2,\ u_3=s_3s_2,\ u_4=s_0s_2,\  u_5=\tau ,\ u_6=\tau s_2,\ u_7=\tau s_3s_2,\ u_8=\tau s_0s_2.$

 According to [D, Theorem 6.4], we have

(a) $ c=\{ u_ix_ku_j^{-1}\,|\, 1\le i,j\le 8, \ k=0,1,2,3,...\}.$

(b) $\Gamma_j=\{ u_ix_ku_j^{-1}\,|\, 1\le i \le 8, \ k=0,1,2,3,...\},\ j=1,2,3,4,5,6,7,8.$

Let $Y=\{1,\ 2, \ ...,\ 7,\ 8\}$ and let ${F_c}$ act on $Y$ trivially. Then $K_{F_c}(Y\times Y)$ is isomorphic to the $8\times 8$ matrix ring $M_8(\Rep {F_c})$, where $\Rep {F_c}$ is the representation ring of   ${F_c}= SL_2(\mathbb C)$.

The main result in this section is the following.

\medskip

{\bf Theorem 4.2.} Let $c$ be the two-sided cell of $W$ (the extended affine Weyl group attached to $Sp_6(\mathbb C)$) containing $s_1s_3$. Then the map
$$\pi: c\to M_8(\Rep F_c ), \quad u_ix_ku_j^{-1}\mapsto V(k)_{ij} $$
induces a ring isomorphism
$$\pi: J_c\to M_8(\Rep F_c ),\quad t_{u_ix_ku_j^{-1}}\mapsto V(k)_{ij},$$
where $V(k)_{ij}$ is the matrix in $M_8(\Rep F_c )$ whose entry at $(p,q)$ is $V(k)$ if $(p,q)=(i,j)$ and is 0 otherwise.

\medskip

{\sl Remark:} The Theorem 4 in [BO] implies that Lusztig's conjecture on the structure of $J_c$ is true. Since under the isomorphism $K_{F_c}(Y\times Y)\simeq M_8(\Rep F_c )$, irreducible ${F_c}$-vector bundles on $Y\times Y$ correspond to $V(k)_{ij}$'s, Theorem 4.2 provides a computable verification for {Lusztig's} conjecture on the structure  of $J_c$.

\medskip

We prove Theorem 4.2 by establishing two lemmas.

\medskip

{\bf Lemma 4.3.} Let $1\le i,j,m,n\le 8$ and $k,l$ be nonnegative integers. Then

(a) $\gamma_{u_ix_ku_j^{-1},u_mx_lu_n^{-1},z}=0\quad \text{if}\ j\ne m\ \text{or}\ x\ne u_i x_pu^{-1}_n\   \text{for some }p;$

(b) $\gamma_{u_ix_ku_j^{-1},u_jx_lu_n^{-1},u_ix_pu_n^{-1}}=\gamma_{x_k,x_l,x_p},\quad\text{for any nonnegative integer } p.$

\medskip

{\bf Proof.} Note that $x_l^{-1}=x_l$. If $ \gamma_{u_ix_ku_j^{-1},u_mx_lu_n^{-1},z}\ne0$, then by 1.4(e) we get $u_ix_ku_j^{-1}\el(u_mx_lu_n^{-1})^{-1}=u_nx_lu_m^{-1}$, $u_ix_ku_j^{-1}\er z,$ and $ u_mx_lu_n^{-1}\el z$. By (b) in Subsection 4.1 we see that the first assertion is true.

Now we prove the second assertion. Let $ *=\{s_1,s_2\},\ \#=\{s_2,s_3\}$ and $\star=\{s_0,s_2\}$. Then

\medskip

(c) $\Gamma_2=\Gamma_1^*,\quad \Gamma_3=\Gamma_1^{\#},\quad \Gamma_4=\Gamma_2^\star,\quad \Gamma_5=\tau \Gamma_1\tau,\ \Gamma_6=\tau \Gamma_2\tau,\ \Gamma_7=\tau \Gamma_3\tau,\ \Gamma_8=\tau \Gamma_4\tau.$

\medskip

Applying   1.4 (j) and 1.4 (k) (repeatedly if necessary)  we get the following identity.
$$\gamma_{u_ix_ku_j^{-1},u_jx_lu_n^{-1},u_ix_pu_n^{-1}}=\gamma_{x_ku_j^{-1},u_jx_l,x_p}.$$
Note that $\tau^2=e$. Again using  1.4 (j) and 1.4 (k)  (repeatedly if necessary) we get the following identity.
$$\gamma_{x_ku_j^{-1},u_jx_l,x_p}=\gamma_{x_k,x_l,x_p}.$$
Part (b) is proved and the lemma is proved. \qed

\medskip

{\bf Lemma 4.4.}  For nonnegative integers {$k,l$,} we have
$$t_{x_k}t_{x_l}=\sum_{0\le p\le\min\{ k,l\}}t_{x_{k+l-2i}}.$$

\medskip

{\bf Proof.} If $k=0$ or $l=0$, the identity above is trivial since $x_0$ is a distinguished involution.

Now assume that $k=1$ and $l\ge 1$.  {Put} $\xi=q^{\frac12}+q^{-\frac12}$.
 Let $H^{<13}$ be the $\mathcal A$-submodule of $H$ spanned by all $C_w$ with $a(w)\ge 3$. By Subsection 1.2 and 1.4 (b) we know that $H^{< 13}$ is a two-sided ideal of $H$. Before continuing, we make a convention: {\sl we shall use the symbol $\Box$ for any element in the two-sided ideal $H^{< 13}$ of $H$. Then $\Box+\Box=\Box$ and $h\Box=\Box$ for any $h\in H$.}

 {First} we have
$$C_{x_1}=C_{\tau s_3s_0s_2s_1s_3}=C_{\tau  s_3}C_{s_0}C_{s_2}C_{s_1}C_{s_3}-C_{\tau s_0s_1s_3},\ \ \text{and}\ C_{\tau s_0s_1s_3}\in H^{<13}.$$

Note that $ {L}(x_k)=\{s_1, s_3\}$. Hence

\begin{align}C_{x_1}C_{x_l}=C_{\tau  s_3}C_{s_0}C_{s_2}C_{s_1}C_{s_3}C_{x_l}+\Box=\xi^2C_{\tau  s_3}C_{s_0}C_{s_2}C_{x_l}+\Box\in \xi^2C_{\tau s_3}C_{s_0}C_{s_2}C_{x_l}+H^{<13}.\end{align}

We compute $C_{\tau s_3}C_{s_0}C_{s_2}C_{x_l}$ step by step.

\medskip

{\bf Step 1.} Compute $C_{s_2}C_{x_l}.$  We have
$$C_{s_2}C_{x_l}=C_{s_2x_l}+\sum\limits_{\substack{z\prec x_l\\ s_2z<z}}\mu(z, x_l)C_z.$$

Note that $L(x_l)=\{s_1,s_3\}$. Assume that $z\prec x_l$ and $s_2z\le z$. If $s_1z\le z$, then $\{s_1,s_2\}\subset L(z)$ and $a(z)\ge a(s_1s_2s_1)=3$. In this case, we have $C_z\in H^{<13}$. If $s_1z\ge z$,  {by 1.1(c)} we must have $s_1z=x_l$. Then $z=\tau s_3s_2x_{l-1}$. This contradicts $s_2z\le z$. Therefore we  have
\begin{align} C_{s_2}C_{x_l}=C_{s_2x_l}+\Box\in C_{s_2x_l}+H^{<13}.\end{align}

\medskip

{\bf Step 2.} Compute $C_{s_0}C_{s_2x_l}.$ We have
$$C_{s_0}C_{s_2x_l}=C_{s_0s_2x_l}+\sum\limits_{\substack{z\prec s_2x_l\\ s_0z<z}}\mu(z, s_2x_l)C_z.$$

Note that $L(s_2x_l)=\{s_2\}$. Assume that $z\prec s_2x_l$ and $s_0z\le z$. If $s_2z\le z$, then $\{s_0,s_2\}\subset {L} (z)$ and $a(z)\ge a(s_0s_2s_0)=3$. In this case, we have $C_z\in H^{<13}$. If $s_2z\ge z$, {by 1.1(c)} we must have $z=x_l$. This contradicts $s_0z\le z$.
Therefore we  have
\begin{align} C_{s_0}C_{s_2x_l}=C_{s_0s_2x_l}+\Box\in C_{s_0s_2x_l}+H^{<13}.\end{align}

\medskip

{\bf Step 3.} Compute $C_{\tau s_3}C_{s_0s_2x_l}.$  We have
$$C_{\tau s_3}C_{s_0s_2x_l}=C_{x_{l+1}}+\sum\limits_{\substack{z\prec s_0s_2x_l\\ s_3z<z}}\mu(z, s_0s_2x_l)C_{\tau z}.$$

Assume that $z\prec s_0s_2x_l$ and $s_3z\le z$. Using 1.4 (b), 1.4 (c) and 1.4 (d), we see that $a(\tau z)\ge 2$ , and if $a(\tau z)=2$ then $x_l\el \tau z\er x_1$. We are only concerned with those $C_{\tau z}$ in above summation with $a(\tau z)=2$. Then $\tau z=x_m$ for some $m<l$ and $L(z)=\{s_0,s_3\}$.  Note that $L( s_0s_2x_l)=\{s_0\}$.  We then have  $\mu(z, s_0s_2x_l)={\tilde\mu({}^\star z, {}^\star(s_0s_2x_l))}=\tilde\mu(s_2z, s_2x_l)=\tilde\mu(^*(s_2z), ^*(s_2x_l))=\tilde\mu(s_1s_2z, x_l)$, where $ \star=\{s_0,s_2\},\ *=\{s_1,s_2\}$.
Since $m<l$, we have $\tilde\mu(s_1s_2z, x_l)=\mu(s_1s_2z, x_l)$. Noting that  $s_3s_1s_2z=s_3s_1s_2\tau x_m\ge s_1s_2\tau x_m$ and $s_3x_l\le x_l$, {by 1.1(c)} we see $s_3s_1s_2z=x_l$, which implies that $\tau z=x_{l-1}$.

In conclusion, if $z\prec s_0s_2x_l$ and $s_3z\le z$, then either $C_z\in H^{<13}$ or $z=\tau x_{l-1}$. Hence we have
\begin{align} C_{\tau s_3}C_{s_0s_2x_l}=C_{x_{l+1}}+C_{x_{l-1}}+\Box.\end{align}

Combining formulas (4)-(7) we get
$$C_{x_1}C_{x_l}=\xi^2(C_{x_{l+1}}+C_{x_{l-1}})+\Box\in \xi^2(C_{x_{l+1}}+C_{x_{l-1}})+H^{<13}.$$
Therefore we have $t_{x_1} t_{x_l}=t_{x_{l+1}}+t_{x_{l-1}}.$

 For $k\ge 2$, since $t_{x_k}=t_{x_1}t_{x_{k-1}}-t_{x_{k-2}}$,  we can use induction on $k$ to prove  the lemma. The argument  is completed. \qed

 \medskip

 {\bf Proof of Theorem 4.2.} Combining Lemmas 4.3 and 4.4 we see that Theorem 4.2 is true.

\section{The based ring of the two-sided cell containing $s_1s_2s_1$}

 {{\bf 5.1.}} In this section we consider the two-sided cell in $W$ containing $s_1s_2s_1$. In [QX] we showed that Lusztig's conjecture on the structure of the based ring of the two-sided cell needs modification. In this section we give a description of the based ring. For consistence, we keep the notations in [QY] for the two-sided cell of $W$ containing $s_1s_2s_1$. In particular, we denote $D$ for the two-sided cell of $W$ containing $s_1s_2s_1$.

 According to [D,Figure I, Theorem 6.4], we have the following result.

 \medskip

 (a) There are  12 left cells in the two-sided cell $D$ and a representative of each left cell in $ D$ are:
$$\begin{array}{lllllllll} &D_{013},&013;&D_2,&0132;&D_{02},&01320;& D_{12},&01321;\\
&&&&&&&&\\
 &D_{3},&01323;&D_{03},&013203; &D_{01},&013201; &D_{13},&013213; \\
&&&&&&&&\\
&D'_{2},&0132032; &\widehat{D'_{2}},&0132132; &D_{1},&01320321; &D_{0},&01321320.   \end{array}$$

\medskip

 The value of $a$-function on $D$ is 3.

 \medskip
 \def\TILDE{\tilde}

 Let $\Gamma$ and $\Gamma'$ be two left cells of $W$. If $\Gamma'=\Gamma^*$ for some $*=\{s,t\}$ (see Subsection 1.3 for definition of $*$-operation), then we write $\Gamma\ \overset{\{s, t\}}{\text{------}}\ \Gamma'$. The following result is easy to verify.

 \medskip

 {\bf Lemma 5.2.} Keep the notations  as above. Then we have

 $$ D_{3}\ \overset{\{s_2, s_3\}}{\text{------}}\ D_{013}\ \overset{\{s_1, s_2\}}{\text{------}}\ D_{2}.$$

$${ \begin{aligned}
D_{0}\ \overset{\{s_0, s_2\}}{\text{------}}\ &\widehat{D'_{2}}\ \overset{\{s_2, s_3\}}{\text{------}}\ D_{12}\ \overset{\{s_0, s_2\}}{\text{------}}\ D_{01}\ \overset{\{s_1, s_2\}}{\text{------}}\ D_{02}\ \overset{\{s_2, s_3\}}{\text{------}}\ D'_{2}\ \overset{\{s_1, s_2\}}{\text{------}}\ D_{1}.\\
&\ | {\scriptstyle{\{s_1, s_2\}}}\qquad\qquad\qquad\qquad\qquad\qquad\qquad\quad\ |{\scriptstyle{\{s_0, s_2\}}}\\
&D_{13}\qquad\qquad\qquad\qquad\qquad\qquad\qquad\qquad\  D_{03}
\end{aligned}}$$

\medskip

\medskip

{\bf 5.3.}  Let $$\begin{array}{lllllll}
u_k=(s_0s_1s_3s_2)^ks_0s_1s_3, &&\\
x_k=(s_1s_2s_3s_0)^ks_1s_2s_1,& \ x'_0=\tau s_2s_0s_1s_2s_1,&\ x'_{k+1}=\tau s_0s_2s_3s_0x_k,  \\
p_1=e,\ & p_2=s_2,\ & p_3=s_3s_2,\\
 p_4=s_1s_2,\ & p_5=s_0s_2,\ &  p_6=s_0s_1s_2,\\
  p_7=s_3s_1s_2,\ & p_8=s_3s_0s_2,\ & p_9=s_2s_3s_1s_2,\\
    p_{10}=s_2s_3s_0s_2,\ & p_{11}=s_0s_2s_3s_1s_2,\ & p_{12}=s_1s_2s_3s_0s_2;\\
 q_4=p_1=e,\ & q_5=\tau,\ & q_6=s_0,\\
  q_7=s_3,\ & q_8=s_3\tau,\ &  q_9=s_2s_3,\\
    q_{10}=s_2s_3\tau,\ & q_{11}=s_0s_2s_3,\ & q_{12}=s_1s_2s_3\tau.\end{array}$$

 According to [D, Theorem 6.4], we have

\medskip

(a) The two-sided cell $D$ consists of the following elements:
$$ p_iu_kp_j^{-1},\   p_i \tau u_kp_j^{-1}, \ q_lx_0q_m^{-1},\ q_lx_0q_6^{-1},\ q_lx_0q_6^{-1}\tau,\ q_lx'_{0}q_m^{-1},$$
where $  1\le i,j\le 12, \ 4\le l\le 12, \ 4\le m\ne 6\le 12,\  k\ge 0.$

\medskip

For convenience, we number the left cells in $D$ as follows:
$$\begin{array}{lllllll}
\Gamma_1=D_{013},& \ \Gamma_2=D_{2},&\ \Gamma_3=D_{3}, & \ \Gamma_4=D_{12}, &\
\Gamma_5=D_{02}, &\ \Gamma_6=D_{01},\\
 \Gamma_7=D_{13}, &\ \Gamma_8=D_{03},&\ \Gamma_9=\widehat{D'_{2}}, &\ \Gamma_{10}=D'_{2}, &\  \Gamma_{11}=D_{0}, &\ \Gamma_{12}=D_{1}.\end{array}$$
Then ({\it loc. cit}) we have

\medskip

(b1) For $j=1,\ 2,\ 3$, the left cell $\Gamma_j$ consists of the following elements:
$$ p_iu_kp_j^{-1},\quad  p_i\tau u_kp_j^{-1},\qquad 1\le i \le 12, \ k\ge 0.$$

\medskip

(b2)  For $j=4,\ 5,\  7,\ 8,\ ...,\ 12$, the left cell $\Gamma_j$ consists of the following elements:
$$p_iu_kp_j^{-1},\quad p_i\tau u_kp_j^{-1},\quad q_lx_0q_j^{-1}, \quad q_lx'_0q_j^{-1},\qquad  { 1\le i\le 12},\ 4\le l\le 12, \ k\ge 0.$$
Note that $p_4u_kp_4^{-1}=x_{k+1},\ p_4\tau u_kp_4^{-1}=x'_{k+1}.$

\medskip

(b3) The left cell $\Gamma_6$ consists of the following {elements:}
$$p_iu_kp_6^{-1},\quad p_i\tau u_kp_6^{-1},\quad q_lx_0q_6^{-1}, \quad q_lx_0q_6^{-1}\tau,\qquad  { 1\le i\le 12},\ 4\le l\le 12, \ k\ge 0.$$

\medskip
\def\Rep{\text{Rep\,}}

{\bf 5.4.} For the two-sided cell $D$, the corresponding unipotent class in $Sp_6(\mathbb C)$ has Jordan block sizes (222). {Maximal}
reductive subgroup of the centralizer of an element in the unipotent class
is ${F_c} = O_3(\mathbb C)=\mathbb Z/2\mathbb Z\times SO_3(\mathbb C)$. Let $\TILDE F_c =\mathbb Z/2\mathbb Z\times SL_2(\mathbb C)$ be the simply connected covering of $F_c$.

Let $Y$ be a set of 12 elements and let $\TILDE F_c$ act on $Y$ trivially. Then $K_{\TILDE F_c }(Y\times Y)$ is isomorphic to the $12\times 12$ matrix ring $M_{12}(\Rep \TILDE F_c )$, where $\Rep \TILDE F_c $ is the representation ring of $\tilde F_c$.

For nonnegative integer $k$, let $V(k)$ be an irreducible representation of $SL_2(\mathbb C)$ with highest weight $k$.  Let $\epsilon$ be the sign representation of $\mathbb Z/2\mathbb Z$. Regarding $V(k)$ and $\epsilon$ as   representations of $\tilde F_c$ naturally, then, up to isomorphism, the irreducible representations of $\tilde F_c$ are $V(k),\ \epsilon V(k),\ k=0,\ 1,\ 2,\ ....$ When $k$ is even, $V(k)$ and $\epsilon V(k)$ are also irreducible representation of $F_c$.

Let $V(k)_{ij}\in M_{12}(\Rep\tilde F_c)$ be the matrix whose entry at $(i,j)$ is $V(k)$  and is 0 elsewhere. Similarly we define $\epsilon V(k)_{ij}$.  The main result in this section is the following theorem.
\medskip

 {\bf Theorem 5.5.}
  There  is a natural injection
\begin{alignat*}{2}
 \pi: c\hookrightarrow& M_{12}(\Rep\ \TILDE F_c ),\\  p_iu_kp_j^{-1}\longmapsto &\begin{cases} {V(2k)_{ij},}&{1\le i, j\le3,}\\
{V(2k+2)_{ij},}&{4\le i, j\le 12,}\\
 {V(2k+1)_{ij},}&{\text{otherwise};}\end{cases}\\
 p_i\tau u_kp_j^{-1}{\longmapsto} &\begin{cases} {\epsilon V(2k)_{ij}, }&{1\le i, j\le3,}\\
 {\epsilon V(2k+2)_{ij},}&{4\le i, j\le 12,}\\
 {\epsilon V(2k+1)_{ij},}&{\text{otherwise};}\end{cases}\\
  y {\longmapsto} & V(0)_{lm}, \quad \text{if $y$ can be obtained from $x_0$ by a}\\
  &\qquad\qquad\quad\text{sequence of left and/or right star operations,}\\
  y {\longmapsto} & \epsilon V(0)_{lm},\quad \text{if $y$ can be obtained from $x'_0$ by a}\\
  &\qquad\qquad\quad\text{sequence of left and/or right star operations,}\\
 \end{alignat*}
 where $y=q_lx_0q_m^{-1}$ or $q_lx'_0q_m^{-1}\ (m\ne 6)$ or $q_lx_0q_6^{-1}$ or $q_lx_0q_6^{-1}\tau$,\ $4\le l,m\le 12$.

The injection $\pi$ induces an injective ring homomorphism
$$\Pi: J_{c}\rightarrow M_{12}(\Rep\tilde F_c)\simeq K_{\tilde F_c}(Y\times Y), \quad t_w\mapsto \pi(w), $$
 where $Y$ is a set of 12 elements with trivial $\TILDE F_c$ action.

\medskip

{\bf Proof:} We need to prove that
\begin{equation}\Pi(t_wt_u)=\pi(w)\cdot \pi(u),\quad \text{for all } w,u\in D.\end{equation}

Since $D$ is the union of all $\Gamma_i,\ 1\le i\le 12$ and $D= D^{-1}$, we know that $D$ is the union of all $\Gamma_i^{-1}\cap\Gamma_j,\ 1\le i,j\le 12.$

Assume that $w\in \Gamma_i^{-1}\cap\Gamma_j$ and $u\in \Gamma_k^{-1}\cap\Gamma_l$. Using 1.4(j), 1.4(k) and Lemma 5.2, we know that  it suffices to prove formula (8) for $w\in \Gamma_i^{-1}\cap\Gamma_j$, $u\in \Gamma_k^{-1}\cap\Gamma_l$, $ i,j,k,l\in\{1, 4\}$. When $j\ne k$, by 1.4 (e) we see that $t_wt_u=0$, hence formula (8) holds in this case. When $i=j=k=l$, according to Theorem 3.1 in [QX], we know that formula (8) holds in this case. To complete the proof of the theorem we need to prove formula (8) for the following cases:

(i) $i=1,\ j=k=1,\ l=4$;

(ii) $i=1,\ j=k=4,\ l=4;$

(iii) $i=1,\ j=k=4,\ l=1;$

(iv) $i=4,\ j=k=1,\ l=1;$

(v) $i=4,\ j=k=1,\ l=4;$

(vi) $i=4,\ j=k=4,\ l=1.$

Keep the notations in the above paragraph. Applying 1.4 (g) and 1.4 (i) we see that to prove formula (8) we only need to prove it for the following two cases: ($\clubsuit$) $w\in \Gamma_4^{-1}\cap\Gamma_1$ and $u\in \Gamma_1^{-1}\cap\Gamma_1$; ($\spadesuit$) $w\in \Gamma_4^{-1}\cap\Gamma_1$ and $u\in \Gamma_1^{-1}\cap\Gamma_4$.

\medskip

{\bf Lemma $\clubsuit$}: We have

(a) $\Gamma_4^{-1}\cap \Gamma_1=\{s_1s_2u_k,\ s_1s_2\tau u_k\,|\, k\ge 0\}$ and $\Gamma_1^{-1}\cap \Gamma_1=\{u_k,\ \tau u_k\,|\, k\ge 0\}$.

\medskip

(b) $\displaystyle t_{s_1s_2u_k}t_{u_l}=t_{s_1s_2\tau u_k}t_{\tau u_l}=\sum_{0\le i\le \min \{2k+1,2l\}}t_{s_1s_2u_{k+l-i}}.$

\medskip

(c) $\displaystyle t_{s_1s_2\tau u_k}t_{u_l}=t_{s_1s_2u_k}t_{\tau u_l}=\sum_{0\le i\le \min \{2k+1,2l\}}t_{s_1s_2\tau u_{k+l-i}}.$


\begin{proof}
Part (a) is obtained from {5.3 (b1) and 5.3 (b2).}

Now we prove (b). Since $u_0$ is a distinguished involution, (b) is true for $l=0$. 

Assume that $l>0$. {First} we will prove \begin{align}t_{s_1s_2u_0}t_{u_l}=t_{s_1s_2u_l}+t_{s_1s_2u_{l-1}},\quad {\text{for any}}\ l>0.\end{align}

\medskip

Let $\xi=q^{\frac12}+q^{-\frac12}$. By a simple computation we get
\begin{align}
&C_{s_1s_2u_0}=(C_{s_1}C_{s_2}-1)C_{s_0s_1s_3}\\
&C_{s_0s_1s_3}C_{u_l}=\xi^3C_{u_l}.\end{align}

Hence
\begin{align} C_{s_1s_2u_0}C_{u_l}=\xi^3(C_{s_1}C_{s_2}-1)C_{u_l}.\end{align}

 Before continuing, we make a convention: {\sl we shall use the symbol $\Box$ for any element in the two-sided ideal $H^{<013}$ of $H$ spanned by all $C_w$ with $a(w)> 3$. Then $\Box+\Box=\Box$ and $h\Box=\Box$ for any $h\in H$.}

\medskip

In  {[QX, }subsection 3.3, Step 1], we have shown the following identity:
\begin{align}C_{s_2}C_{u_l}=C_{s_2u_l}+\Box\in C_{s_2u_l}+H^{<013}.\end{align}

\medskip

Now we compute $C_{s_1}C_{s_2u_l}$. By formula (1) in 1.1 (e),  we have $$C_{s_1}C_{s_2u_l }=C_{s_1s_2u_l }+\sum\limits_{\substack{y\prec s_2u_l\\ s_1y<y}}\mu(y, s_2u_l )C_y.$$
Note that $ {L} (s_2u_l )=\{s_2\}$. {First} we have $\mu(u_l,s_2u_l)=1$ and $s_1u_l<u_l$. By 1.4 (c) and 1.4 (d), if $C_y$ appears in the above summation with nonzero coefficient, $y\ne u_l$ and $C_y\not\in H^{<013}$, then $y\in\Gamma_4^{-1}\cap \Gamma_1$.

Assume $y\prec s_2u_l$, $s_1y<y$ and $y\in\Gamma_4^{-1}\cap \Gamma_1$. {  By (a) we must have $y{=s_1s_2u_k}=s_1s_2(s_0s_1s_3s_2)^ks_0s_1s_3$ for some nonnegative integer $k\le l-1$.} Since $s_2s_0y\ge s_0y\ge y$, by 1.3 (d) we get $\mu(y,s_2u_l)=\mu(s_0y,u_l)$. Now $s_3u_l\le u_l$ and $s_3s_0y\ge s_0y$, by 1.1(c) we get $s_3s_0y=u_l$. Hence $y=s_1s_2u_{l-1}$.

We have shown
  \begin{align}
  C_{s_1}C_{s_2u_l}=C_{s_1s_2u_l}+C_{u_l}+C_{s_1s_2u_{l-1}}+\Box.\end{align}

Combining formulas (12) (13) and (14), we get formula (9).

{Recall the following formula} in [QX, 3.3]:
\begin{align}t_{u_k}t_{u_l}=\sum_{0\le i\le \min  \{2k,2l\}}t_{u_{k+l-i}}.
\end{align}
Now we employ formulas (9) and (15) to prove the identity in (b). We use induction on $k$. When $k=0$, it is just formula (9). Assume that the formula in (b) is true for nonnegative integer less than $k$. We have
$$\begin{aligned}t_{s_1s_2u_k}t_{u_l}=&(t_{s_1s_2u_0}t_{u_k}-t_{s_1s_2u_{k-1}})t_{u_l}\\=&t_{s_1s_2u_0}\cdot\sum_{0\le i\le \min \{2k,2l\}}t_{u_{k+l-i}}-t_{s_1s_2u_{k-1}}t_{u_l}\\
=&\sum_{0\le i\le \min \{2k,2l\}}t_{s_1s_2u_{k+l-i}}+\sum_{0\le i\le \min \{2k,2l\}}t_{s_1s_2u_{k+l-i-1}}-\sum_{0\le j\le \min \{2k-1,2l\}}t_{s_1s_2u_{k+l-1-j}}\\
=&\sum_{0\le i\le \min \{2k,2l\}}t_{s_1s_2u_{k+l-i}}+\sum_{1\le i\le \min \{2k+1,2l+1\}}t_{s_1s_2u_{k+l-i}}-\sum_{1\le j\le \min \{2k,2l+1\}}t_{s_1s_2u_{k+l-j}}\\
=&\sum_{0\le i\le \min \{2k+1,2l\}}t_{s_1s_2u_{k+l-i}}
\end{aligned}$$

Since $\tau u_k=u_k\tau$, by 1.4 (j) we have $t_{s_1s_2\tau u_k}t_{\tau u_l}=t_{s_1s_2u_k}t_{u_l}$. Part (b) is proved.

Since $\tau u_k=u_k\tau$ and $\tau u_{k+l-i}=u_{k+l-i}\tau$,  using 1.4 (j) we see that  Part (c) follows from Part (b).

The proof is completed.
\end{proof}

\medskip

{\bf Lemma $\spadesuit$.} (a) For $k\ge 0$, we have $s_1s_2u_ks_2s_1=x_{k+1}$ and $s_1s_2\tau u_ks_2s_1=x'_{k+1}$. Moreover, $\Gamma_4^{-1}\cap \Gamma_4=\{x_k,\ x'_k\,|\, k\ge 0\}$.

 For nonnegative integers $k,l$ we have

(b) $\displaystyle t_{s_1s_2u_k}t_{u_ls_2s_1}=t_{s_1s_2\tau u_k}t_{u_l\tau s_2s_1}=\sum_{0\le i\le \min \{2k+1,2l+1\}}t_{x_{k+l+1-i}}.$

\medskip

(c) $\displaystyle t_{s_1s_2\tau u_k}t_{u_ls_2s_1}=t_{s_1s_2  u_k}t_{u_l\tau s_2s_1}=\sum_{0\le i\le \min \{2k+1,2l+1\}}t_{x'_{k+l+1-i}}.$

\medskip

\begin{proof}
Part (a) follows from the discussion in Subsection 5.3.

Now we prove Part (b). {First} we  prove \begin{align}t_{s_1s_2u_0}t_{u_ls_2s_1}=t_{x_{l+1}}+t_{x_l}.\end{align}

In [QX, Subsection 4.2], it is shown $t_{s_1s_2u_0}t_{u_0s_2s_1}=t_{x_1}+t_{x_0}$. Now assume that $l\ge 1$. As before,
 $\xi=q^{\frac12}+q^{-\frac12}$.  Since $C_{s_0s_1s_3}C_{u_ls_2s_1}=\xi^3C_{u_ls_2s_1}$, using formula (10) we get
\begin{align}C_{s_1s_2u_0}C_{u_ls_2s_1}=\xi^3(C_{s_1}C_{s_2}-1)C_{u_ls_2s_1}.\end{align}

\medskip

We compute the right hand side of equality (17) step by step. As in the proof of Lemma $\clubsuit$, {\sl we shall use the symbol $\Box$ for any element in the two-sided ideal  $H^{<013}$ of $H$ spanned by all $C_w$ with $a(w)> 3$.}

\medskip

{\bf Step 1:} Compute $C_{s_2}C_{u_ls_2s_1}$.

By formula (1) in 1.1 (e), we have $C_{s_2}C_{u_ls_2s_1}=C_{s_2u_ls_2s_1}+\sum\limits_{\substack{y\prec u_ls_2s_1\\ s_2y<y}}\mu(y, u_ls_2s_1)C_y.$
Note that $ {L} (u_ls_2s_1)=\{s_0, s_1, s_3\}$.

Assume $y\prec u_ls_2s_1$ and $s_2y<y$.
\begin{itemize}
\item If $s_0y>y$, then by 1.1(c) we get $s_0y=u_ls_2s_1$. This contradicts the assumption. So $s_0y>y$ would not occur.
\item If $s_1y>y$, then by 1.1(c) we get $s_1y=u_ls_2s_1$. This contradicts the assumption. So $s_1y>y$ would not occur.
\item If $s_0y<y, s_1y<y$ and $s_2y<y$, then $a(y)\ge a(w_{012})=6$. So $C_y\in H^{<013}.$
\end{itemize}

Therefore,
\begin{align} C_{s_2}C_{u_ls_2s_1}=C_{s_2u_ls_2s_1}+\Box.\end{align}

\medskip

{\bf Step 2:}
Similar to the proof for formula (14), we have
  \begin{align}C_{s_1}C_{s_2u_ls_2s_1}=C_{s_1s_2u_ls_2s_1}+C_{u_ls_2s_1}+C_{s_1s_2u_{l-1}s_2s_1}+\Box.\end{align}

Note $s_1s_2u_ls_2s_1=x_{l+1}$.
Combining formulas (17)-(19), we get (16).

Now we can  prove part (b) using induction on $k$. By 1.4 (j), we know $t_{s_1s_2u_k}t_{u_ls_2s_1}=t_{s_1s_2\tau u_k}t_{u_l\tau s_2s_1}$. Thus for  $k=0$,  Part (b) is equivalent to formula (16), which is true. Now assume that $k\ge 1$ and Part (b) is true for $k-1$. Using   Lemma $\clubsuit$ and 1.4(i), induction hypothesis and formula (16),  we get
$$\begin{aligned}
t_{s_1s_2u_k}t_{u_ls_2s_1}=&(t_{s_1s_2u_0}t_{u_k}-t_{s_1s_2u_{k-1}})t_{u_ls_2s_1}\\=&{t_{s_1s_2u_0}
\cdot}\sum_{0\le i\le \min \{2l+1,2k\}}t_{u_{k+l-i}s_2s_1}-t_{s_1s_2u_{k-1}}t_{u_ls_2s_1}\\
=&\sum_{0\le i\le \min \{2l+1,2k\}}(t_{x_{k+l+1-i}}+t_{x_{k+l-i}})-\sum_{0\le i\le \min \{2l+1,2k-1\}}t_{x_{k+l-i}}\\
=&\sum_{0\le i\le \min \{2l+1,2k+1\}}t_{x_{k+l+1-i}}.\end{aligned}$$
This completes the proof for Part (b).

Proof for Part (c) is similar. {First}, it is easy to check that
$$C_{s_0s_2u_0}C_{u_0s_2s_1}=\xi^3(C_{\tau x'_1}+C_{\tau x'_0}),$$
which implies
$$t_{s_1s_2\tau u_0}t_{u_0s_2s_1}=t_{x'_1}+t_{x'_0}.$$
Further, we prove that
$$t_{s_1s_2\tau u_0}t_{u_ls_2s_1}=t_{x'_{1+1}}+t_{x'_l}.$$
Then using induction on $k$, as the proof for Part (b), we prove Part (c). The proof for Lemma $\spadesuit$ is completed.
\end{proof}

We have completed the proof for Theorem 5.5.

\medskip

{\bf 5.6.} Motivated by Theorem 5.5 and the discussion of the cocenter of $J$ in [BDD, Section 5] and some other evidences, we suggest a modification of Lusztig's conjecture on the structure of $J_c$ which is stated for any connected reductive groups over $\mathbb C$.

Let $W$ be the extended affine Weyl group attached to a connected reductive group over $\mathbb C$ (see Subsection 1.1) and let $c$ be a two-sided cell of $W$. Let $F_c$ be a maximal reductive subgroup of the centralizer of an element in the corresponding unipotent class of $G$. Then there should exist a reductive group $\tilde F_c$ with the following properties:

(i) The reductive group $\tilde F_c$ is a simply connected covering of $F_c$. That is, the identity component $\tilde F_c^\circ$ has  simply connected derived group, and there is a natural surjective homomorphism $\tilde F_c\to F_c$ with finite kernel. In particular, if $  F_c^\circ$ has  simply connected derived group, then $\tilde F_c=F_c$.

(ii) There exists a finite set $Y$ with an algebraic action of $\tilde F_c$ and an injection
$$\pi: c\hookrightarrow\text{isomorphism classes of irreducible $\tilde F_c$-vector bundles on}\ Y\times Y.$$
  such that

  (iii)  The injection $\pi$ induces a ring injection
  $$\Pi: J_c\to K_{\tilde F_c}(Y\times Y),\ \ t_x{\mapsto} \pi(x).$$

  (iv) $\pi(x^{-1})_{(a,b)}=\pi(x)_{(b,a)}^*$ is the dual representation of $\pi(x)_{(b,a)}.$

  (v) $K_{\tilde F_c}(Y\times Y)$ is a finitely generated left (and right as well) $\Pi(J_c)$-module.

\medskip

 It seems natural that the $F_c$-set $\mathbf B_e$  defined in a recent paper (see [L4]) would have an $\tilde F_c$-action compatible with the $F_c$-action and then $\mathbf B_e$ should be a good candidate for the set $Y$ above.

\bigskip

\noindent{\bf Acknowledgement:} Part of the work was done during
YQ's visit to the Academy of Mathematics and Systems Science, Chinese Academy of Sciences. YQ
is very grateful to the AMSS for hospitality and for
financial supports.



\begin{thebibliography}{99}

\bibitem[B]{B} R. Bezrukavnikov, {\sl On tensor categories attached to cells in affine Weyl groups,} In "Representation Theory of Algebraic Groups and Quantum Groups", Adv. Stud. Pure Math., 40, Math. Soc. Japan, Tokyo, 2004, pp. 69-90.

\bibitem[BDD]{BDD} R. Bezrukavnikov,  S. Dawydiak, G. Dobrovolska, {\sl On the structure of the affine asymptotic Hecke algebras}, arXiv:2110.15903.

\bibitem[BO]{BO} R. Bezrukavnikov and V. Ostrik, {\sl On tensor categories attached to cells in affine Weyl groups II,} In "Representation Theory of Algebraic Groups and Quantum Groups", Adv. Stud. Pure Math., 40, Math. Soc. Japan, Tokyo, 2004,  pp.101-119.

\bibitem[DLP]{DLP} C. De Concini,  G. Lusztig,   C. Procesi,{\sl  Homology of the zero-set of a nilpotent vector field on a flag manifold}, J. Amer. Math. Soc. 1 (1988), 15-34.


\bibitem[D]{D} J. Du, {\sl The decomposition into cells of the affine Weyl group of type $\tilde{B_3}$,} Communications in Algebra, 16 (1988), no.7, 1383--1409.


\bibitem[KL]{KL} D. Kazhdan and G. Lusztig, {\sl Representations of Coxeter
groups and Hecke algebras,} Invent. Math. 53 (1979), 165-184.

\bibitem[L1]{L1} G. Lusztig, {\sl Cells in affine Weyl groups,} in
``Algebraic groups and related topics", Advanced Studies in Pure
Math., vol. 6, Kinokunia and North Holland, 1985, pp. 255-287.

\bibitem[L2]{L2} G. Lusztig, {\sl Cells in affine Weyl groups, II,} J. Alg.
109 (1987), 536-548.

\bibitem[L3]{L3} G. Lusztig, {\sl Cells in affine Weyl groups, IV,} Journal of The Faculty of Science, 36 (1989), no.2, 297-328.

    \bibitem[L4]{L4} G. Lusztig, {\sl Discretization of Speinger fibers, } arXiv:1712.07530v3, 2021.

 \bibitem[QX]{QX} Yannan Qiu and Nanhua Xi, {\sl The based ring of two-sided cells in an affine Weyl group of type $\tilde B_3$, I.} Sci. China Math., to appear. arxiv: 2107.08983.




\bibitem[X1]{X1} N. Xi, {\sl Representations of Affine Hecke Algebras,} volume 1587, Springer Lecture Notes in Math., 1994.

\bibitem[X2]{X2} N. Xi, {\sl The based ring of two-sided cells of affine Weyl groups of type ${\tilde{A}_{n-1}}$,} volume 749, American Mathematical Soc., 2002.






\end{thebibliography}
\end{document}